# Approximation of Sums of Primes.

**Nikos Bagis**

Department of Informatics

Aristotle University of Thessaloniki

54124 Thessaloniki, Greece

bagkis@hotmail.com

**Abstract.**

In this work we consider sums of primes that converging very slow. We set as a base, a reformulation of analytic prime number theorem and we use the values of Riemann Zeta function for the approximation. We also give the truncation error of these approximations.

## §1. Finite Sums and Products

**Theorem 1.** Let $f$ be analytic in $A \supseteq (-1,1]$ with $f(0) = 0$. Let also $a_n$ is an arbitrary sequence such that $|a_n| < 1$, for $n = 1, 2, 3, \ldots$. Then we have

$$\sum_{k \leq x} f(a_k) = \sum_{n=1}^{\infty} \frac{\log\left(\prod_{k \leq x}(1 - a_k^n)\right)}{n} \sum_{d|n} \frac{f^{(d)}(0)}{\Gamma(d)} \mu\left(\frac{n}{d}\right) \quad :(1)$$

Where is a positive integer and $\mu$ is the Moebius Mu function defined by

$$\mu(k) = \begin{cases} 1, & \text{if } k = 1 \\ (-1)^r, & \text{if } k = \text{the product of } r \text{ distinct primes} \\ 0, & \text{otherwize} \end{cases}$$

**Proof.**

$$\sum_{n=1}^{\infty} \frac{\log\left(\prod_{k \leq x}(1 - a_k^n)\right)}{n} c_n = \sum_{n=1}^{\infty} \frac{c_n}{n} \sum_{k=1}^{x} \log(1 - a_k^n) = \sum_{n=1}^{\infty} \frac{c_n}{n} \sum_{k=1}^{x} \sum_{m=1}^{\infty} \frac{1}{m} a_k^{mn} = \sum_{k=1}^{x} \sum_{n=1}^{\infty} \sum_{m=1}^{\infty} \frac{c_n}{nm} a_k^{mn} =$$

$$\sum_{k=1}^{x} \sum_{n=1}^{\infty} \sum_{m=1}^{\infty} \frac{c_n}{nm} a_k^{mn} \stackrel{nm=r}{=} \sum_{k=1}^{x} \sum_{r=1}^{\infty} \frac{1}{r} a_k^r \sum_{d|r} c_d \quad :(2)$$

Let now $\dfrac{1}{r}\sum_{d|r} c_d = \dfrac{f^{(r)}(0)}{r!}$, from the Moebius Inversion Theorem (see **[A]** page 54),

we get that $c_r = \sum_{d|r} \dfrac{f^{(d)}(0)}{\Gamma(d)} \mu\left(\dfrac{r}{d}\right)$. Hence equation (2) becomes:

$$\sum_{n=1}^{\infty} \dfrac{\log\left(\prod_{k\leq x}(1-a_k^n)\right)}{n} \sum_{d|n} \dfrac{f^{(d)}(0)}{\Gamma(d)} \mu\left(\dfrac{n}{d}\right) = \sum_{k=1}^{x} \sum_{r=1}^{\infty} \dfrac{f^{(r)}(0)}{r!}(a_k)^r = \sum_{k=1}^{x} f(a_k). \quad \square$$

## §2. The extension of Euler's Theorem

**Theorem 2.** If $f$ is analytic in $A \supseteq (-1,1]$ and $f(0)=0$, then we have

$$\sum_{p} f\left(\dfrac{1}{p^s}\right) = \sum_{n=1}^{\infty} \dfrac{\log(\zeta(ns))}{n} \sum_{d|n} \dfrac{f^{(d)}(0)}{\Gamma(d)} \mu\left(\dfrac{n}{d}\right) \quad :(3)$$

Where the first sum is over all primes, and $\zeta$ is the Riemann's Zeta function:
$$\zeta(s) = \sum_{n=1}^{\infty} \dfrac{1}{n^s}, \; s>1$$

**Proof.** Set $a_n = 1/p_n^s$, where $p_n$ is the $n$-th prime, $s>1$ and $x=\infty$. Then from Euler's Theorem $1/\zeta(s) = \prod_{p-prime}\left(1-\dfrac{1}{p^s}\right)$, $s>1$ the result follows easily from Theorem 1. $\square$

**Proposition 1.** The Truncation Error of (3) is:

$$E(f,M,s) := \left| \sum_{n=M+1}^{\infty} \sum_{d|n} \dfrac{f^{(d)}(0)}{\Gamma(d)} \mu\left(\dfrac{n}{d}\right) \dfrac{\log(\zeta(ns))}{n} \right|,$$

Then
$$E(f,M,s) \leq \dfrac{2^{s+1}(2^s+1)(s+1)}{\pi(2^s-1)^3} \int_0^{2\pi} |f(e^{it})| dt \dfrac{M^2}{(sM+s-1)2^{sM}}$$

**Proof.** Let $E = \left| \sum_{n=M+1}^{\infty} \sum_{d|n} \dfrac{f^{(d)}(0)}{\Gamma(d)} \mu\left(\dfrac{n}{d}\right) \dfrac{\log(\zeta(ns))}{n} \right|$

First from the Poisson formula we have:
$\dfrac{f^{(k)}(0)}{k!} = \dfrac{1}{2\pi} \int_0^{2\pi} f(e^{it}) e^{-itk} dt$, thus

$$\left|\frac{f^{(k)}(0)}{\Gamma(k)}\right| \le \frac{k}{2\pi}\int_0^{2\pi}\left|f(e^{it})\right|dt$$

$$\left|\sum_{d|n}\frac{f^{(d)}(0)}{\Gamma(d)}\mu\left(\frac{n}{d}\right)\right| \le \sum_{k=1}^{n}\left|\frac{f^{(k)}(0)}{\Gamma(k)}\right| \le \frac{1}{2\pi}\int_0^{2\pi}\left|f(e^{it})\right|dt\sum_{k=1}^{n}k$$

and the error becomes

$$E \le \frac{1}{4\pi}\int_0^{2\pi}\left|f(e^{it})\right|dt\left|\sum_{n=M+1}^{\infty}(n+1)\log(\zeta(sn))\right| \le$$

$$\frac{1}{4\pi}\int_0^{2\pi}\left|f(e^{it})\right|dt\left|\sum_{n=M+1}^{\infty}(n+1)\left(\frac{1}{2^{sn}}+\frac{1}{3^{sn}}+\frac{1}{4^{sn}}+\ldots\right)\right| \le$$

(From the inequality $\log(x) \le x-1$).

$$\frac{1}{4\pi}\int_0^{2\pi}\left|f(e^{it})\right|dt\left|\sum_{n=M+1}^{\infty}(n+1)\left(\frac{1}{2^{sn}}+\frac{1}{3^{sn}}+\frac{1}{4^{sn}}+\ldots\right)\right| \quad :(4)$$

We try to estimate the series $A(x) := \frac{1}{2^x}+\frac{1}{3^x}+\frac{1}{4^x}+\ldots$, $x = sn$

$$A(x) := \frac{1}{2^x}+\frac{1}{3^x}+\frac{1}{4^x}+\ldots \le \frac{1}{2^x}+\int_2^{\infty}\frac{1}{t^x}dt =$$

$$\frac{1}{2^{sn}}+\frac{2}{(sn-1)2^{sn}} = \frac{1}{2^{sn}}\frac{sn+1}{sn-1}$$

and (4) can be estimated by:

$$\frac{1}{4\pi}\int_0^{2\pi}\left|f(e^{it})\right|dt\left|\sum_{n=M+1}^{\infty}(n+1)\left(\frac{1}{2^{sn}}+\frac{1}{3^{sn}}+\frac{1}{4^{sn}}+\ldots\right)\right| \le$$

$$\frac{1}{4\pi}\int_0^{2\pi}\left|f(e^{it})\right|dt\left|\sum_{n=M+1}^{\infty}(n+1)\frac{1}{2^{sn}}\frac{sn+1}{sn-1}\right|$$

Observe that: $n+1 \le 2n$, $\frac{1}{sn-1} \le \frac{1}{s(M+1)-1}$,

thus

$$\frac{1}{4\pi}\int_0^{2\pi}\left|f(e^{it})\right|dt\left|\sum_{n=M+1}^{\infty}(n+1)\frac{1}{2^{sn}}\frac{sn+1}{sn-1}\right| \le \frac{2}{4\pi(sM+s-1)}\int_0^{2\pi}\left|f(e^{it})\right|dt\sum_{n=M+1}^{\infty}\frac{n(sn+1)}{2^{sn}} \le$$

$$\frac{(s+1)}{2\pi(sM+s-1)}\int_0^{2\pi}\left|f(e^{it})\right|dt\sum_{n=M+1}^{\infty}\frac{n^2}{2^{sn}} =$$

$$\frac{(s+1)}{2\pi(sM+s-1)}\int_0^{2\pi}\left|f(e^{it})\right|dt\frac{1}{2^{Ms}}\sum_{n=1}^{\infty}\frac{(n+M)^2}{2^{sn}} =$$

$$\frac{(s+1)}{2\pi(sM+s-1)}\int_0^{2\pi}\left|f(e^{it})\right|dt\frac{1}{2^{Ms}}\sum_{n=1}^{\infty}\frac{n^2+2nM+M^2}{2^{sn}} \le$$

$$\frac{(s+1)}{2\pi(sM+s-1)}\int_0^{2\pi}|f(e^{it})|dt\frac{1}{2^{M\cdot s}}\sum_{n=1}^{\infty}\frac{4M^2n^2}{2^{sn}} \leq$$

$$\frac{(s+1)2M^2}{\pi(sM+s-1)}\int_0^{2\pi}|f(e^{it})|dt\frac{1}{2^{Ms}}\sum_{n=1}^{\infty}\frac{n^2}{2^{sn}} =$$

$$\frac{2^{s+1}(2^s+1)(s+1)}{\pi(2^s-1)^3}\int_0^{2\pi}|f(e^{it})|dt\frac{M^2}{(sM+s-1)2^{sM}}.$$

And we get the estimate. □

## §3. The Euler Totient constant

When $n$ is a positive integer, Euler's Totient function, $\varphi(n)$, is defined to be the number of positive integers not greater than $n$ and relatively prime to $n$. Interesting constants emerge if we consider the sum of reciprocals of $\varphi(n)$. Landau proved that:

$$\sum_{n=1}^{N}\frac{1}{\phi(n)} = A\log(N) + B + O\left(\frac{\log(N)}{N}\right),$$

where:

$$A = \sum_{k=1}^{\infty}\frac{\mu(k)^2}{k\phi(k)} = \frac{\zeta(2)\zeta(3)}{\zeta(6)} = \frac{315}{2\pi^4}\zeta(3) = 1.9435964368...$$

and:

$$B = \frac{315}{2\pi^4}\zeta(3) - \sum_{k=1}^{\infty}\frac{\mu(k)^2\log(k)}{k\phi(k)} = -0.06057... \ .$$

In the above formulas, $\zeta(x)$ is Riemann's zeta function, $\zeta(3)$ is Apery`s

One proof of the above asymptotic result can be found in Koninck & Ivic.

An alternative expression for the constant $B$ is:

$$B = \frac{315}{2\pi^4}\zeta(3)\left(\gamma - \sum_p\frac{\log(p)}{p^2-p+1}\right) = -0.06057...$$

We will see how to estimate the value of some prime series, such as:

$$\sum_p\frac{\log(p)}{p^2-p+1} \quad :(a)$$

For this we set

$$f_1(x) = \frac{1}{\sqrt{3}}\arctan\left(\frac{-1+2x}{\sqrt{3}}\right) + \frac{1}{2}\log(1-x+x^2) \quad :(5),$$

then

$$B = \frac{315}{2\pi^4} \zeta(3) \left( \gamma + \sum_{n=2}^{\infty} \frac{\zeta'(n)}{\zeta(n)} \sum_{d|n} \frac{f_1^{(d)}(0)}{\Gamma(d)} \mu\left(\frac{n}{d}\right) \right).$$

(We can calculate $\sum_p \frac{\log(p)}{p^2 - p + 1}$, with differentiation of (3), with respect to $s$.

$$\sum_p f_1'\left(\frac{1}{p^s}\right) \frac{\log(p)}{p^s} = \sum_p \frac{\log(p)}{p^2 - p + 1} \text{ and set } s = 1)$$

**Proposition 2.** If we set

$$E_1(f_1, M) := \left| \sum_{n=M+1}^{\infty} \frac{\zeta'(n)}{\zeta(n)} \sum_{d|n} \frac{f_1^{(d)}(0)}{\Gamma(d)} \mu\left(\frac{n}{d}\right) \right|, M \gg 1.$$

then

$$E(f_1, M) \leq (2M + 4)\zeta'(M + 1),$$

where $\zeta'$ is the first derivative of the Riemann's Zeta function.

**Proof.** If $a(k) := \frac{f_1^{(k)}(0)}{k!} k$, then

$$a(k) = \begin{cases} 0, & \text{if } k \equiv 1 \pmod{6} \\ 1, & \text{if } k \equiv 2 \pmod{6} \\ 1, & \text{if } k \equiv 3 \pmod{6} \\ 0, & \text{if } k \equiv 4 \pmod{6} \\ -1, & \text{if } k \equiv 5 \pmod{6} \\ -1, & \text{if } k \equiv 0 \pmod{6} \end{cases}, f(0)=0$$

(Note also that $a(k) = \frac{2}{\sqrt{3}} \sin\left(\frac{(k-1)\pi}{3}\right)$, for $k = 1, 2, 3,\ldots$)

Thus:

$$\left| \sum_{d|n} \frac{f_1^{(d)}(0)}{\Gamma(d)} \mu\left(\frac{n}{d}\right) \right| < \sum_{d|n} |\mu(d)| \leq d(n) \leq n,$$

where $d(n) = \sum_{d|n} 1$.

$$E(f, M) := \left| \sum_{n=M+1}^{\infty} \frac{\zeta'(n)}{\zeta(n)} \sum_{d|n} \frac{f_1^{(d)}(0)}{\Gamma(d)} \mu\left(\frac{n}{d}\right) \right|, M \gg 1$$

Thus

$$E(f, M) \leq \sum_{n=M+1}^{\infty} \left|\frac{\zeta'(n)}{\zeta(n)}\right| \sum_{d|n} 1 \leq \sum_{n=M+1}^{\infty} |\zeta'(n)| \sum_{k=1}^{n} 1 =$$

$$\sum_{n=M+1}^{\infty} \left( \frac{\log 2}{2^n} + \frac{\log 3}{3^n} + \ldots \right) n = \sum_{n=M+1}^{\infty} \sum_{k=2}^{\infty} \frac{\log(k)}{k^n} n = \sum_{k=2}^{\infty} \frac{\log(k)}{k} \sum_{n=M}^{\infty} \frac{n+1}{k^n} =$$

$$\sum_{k=2}^{\infty}\frac{\log(k)}{k}\sum_{n=M}^{\infty}\frac{n+1}{k^n} \le \sum_{k=2}^{\infty}\frac{\log(k)}{k^M}\left(\frac{k-M+kM}{(k-1)^2}\right) \le \sum_{k=2}^{\infty}\frac{\log(k)}{k^M}\left(\frac{4}{k}+\frac{2M}{k}\right)=$$
$$(2M+4)\zeta'(M+1). \quad \square$$

**Note:** As you can see the main interest is not in the values of $f$ but in the values of $\zeta(x)$, $\zeta'(x)$, where $x>1$.

### §4. Regularized Products over all Primes

Given an increasing sequence $0 < \lambda_1 \le \lambda_2 \le \lambda_3 \le ...$ one defines the regularized (or zeta-regularized) infinite product as:
$$\coprod_n \lambda_n = \exp(-\zeta'_\lambda(0)),$$
where $\zeta_\lambda$ is the zeta function associated to the sequence $(\lambda_n)$,
$$\zeta_\lambda(s) = \sum_{n=1}^{\infty}\frac{1}{\lambda_n^s}$$

From Theorem 2 we have $\sum_p f\left(\frac{1}{p^s}\right) = \sum_{n=1}^{\infty}\frac{\log(\zeta(ns))}{n}\sum_{d|n}\frac{f^{(d)}(0)}{\Gamma(d)}\mu\left(\frac{n}{d}\right)$.

Differentiating once the above relation and taking $s = 0$ we get:

**Proposition 3.**
$$\coprod_p p = 4\pi^2$$
is valid

**Proof.** The proof can be found also in **[G,M]** . $\square$

Now if we get the derivative of the main theorem two times and set $s = 0$ we have
$$\sum_{p-prime} f'\left(\frac{1}{p^s}\right)\frac{1}{p^s}(-\log(p)) = \sum_{n=1}^{\infty}\frac{\zeta'(sn)}{\zeta(sn)}\sum_{d|n}\frac{f^{(d)}(0)}{d!}d\mu\left(\frac{n}{d}\right) \Rightarrow$$
$$\sum_{p-prime} f''\left(\frac{1}{p^s}\right)\frac{1}{p^{2s}}(-\log(p))^2 + \sum_{p-prime} f'\left(\frac{1}{p^s}\right)\frac{1}{p^s}(-\log(p))^2 =$$
$$= \sum_{n=1}^{\infty}\frac{\zeta''(sn)\zeta(sn)-(\zeta'(sn))^2}{\zeta(sn)^2}n\sum_{d|n}\frac{f^{(d)}(0)}{d!}d\mu\left(\frac{n}{d}\right)$$

Setting $s = 0$ we have
$$\sum_{p-prime} f''(1)(\log(p))^2 + \sum_{p-prime} f'(1)(\log(p))^2 =$$
$$= \sum_{n=1}^{\infty}\frac{\zeta''(0)\zeta(0)-(\zeta'(0))^2}{\zeta(0)^2}n\sum_{d|n}\frac{f^{(d)}(0)}{d!}d\mu\left(\frac{n}{d}\right) =$$
$$(-2\zeta''(0)-\log(2\pi)^2)\sum_{n=1}^{\infty}n\sum_{d|n}\frac{f^{(d)}(0)}{d!}d\mu\left(\frac{n}{d}\right).$$

But

$$\sum_{n=1}^{\infty} n \sum_{d|n} \frac{f^{(d)}(0)}{d!} d\mu\left(\frac{n}{d}\right) = \lim_{h \to -1} \left( \left( \sum_{n=1}^{\infty} \frac{\frac{f^{(n)}(0)}{n!} n}{n^h} \right) \left( \sum_{n=1}^{\infty} \frac{\mu(n)}{n^h} \right) \right) = \frac{f'(1) + f''(1)}{\zeta(-1)}$$

Hence we get the next

**Proposition 4.** The next regularized product is valid:
$$\coprod_p p^{\log(p)} = e^{2\zeta''(0) + 12(\log(2\pi))^2}$$

### §5. Finite Sums of Primes

**Theorem 3.** Let $x$ be an integer greater than 1, $s$ is real: $s > 1$ and $g$ analytic function in $A \supseteq (-1,1]$, such that $\frac{g^{(a)}(0)}{a!}$ are bounded above (see **[V]**). Then

$$\sum_{n \leq x} \frac{\mu(n)}{n} \sum_{k=1}^{\infty} \sum_{l=1}^{\infty} \frac{1}{l} g\left(\frac{1}{p_k^{snl}}\right) = \sum_{k=1}^{\infty} g\left(\frac{1}{p_k^s}\right) + O\left(x \cdot g\left(\frac{1}{2^{s(x+1)}}\right)\right) \quad :(6)$$

**Proof.** Set $f(x) = x$, $s > 1$ in relation of Theorem 2 then

$$\sum_{k=1}^{\infty} \frac{1}{p_k^s} = -\sum_{n=1}^{\infty} \frac{\mu(n)}{n} \log \prod_{k=1}^{\infty} \left(1 - \frac{1}{p_k^{sn}}\right)$$

$$= -\sum_{n=1}^{M} \frac{\mu(n)}{n} \log \left( \prod_{k=1}^{\infty} \left(1 - \frac{1}{p_k^{sn}}\right) \right) - \sum_{n=M+1}^{\infty} \frac{\mu(n)}{n} \log \left( \prod_{k=1}^{\infty} \left(1 - \frac{1}{p_k^{sn}}\right) \right) =$$

$$= -\sum_{n=1}^{M} \frac{\mu(n)}{n} \sum_{k=1}^{\infty} \log\left(1 - \frac{1}{p_k^{sn}}\right) + \sum_{n=M+1}^{\infty} \frac{\mu(n)}{n} \log(\zeta(sn))$$

$$= \sum_{n=1}^{M} \frac{\mu(n)}{n} \sum_{k=1}^{\infty} \sum_{l=1}^{\infty} \frac{1}{l} \frac{1}{p_k^{snl}} + \sum_{n=M+1}^{\infty} \frac{\mu(n)}{n} \log(\zeta(sn))$$

thus

$$\sum_{k=1}^{\infty} \frac{1}{p_k^s} - \sum_{n=1}^{M} \frac{\mu(n)}{n} \sum_{k=1}^{\infty} \sum_{l=1}^{\infty} \frac{1}{l} \frac{1}{p_k^{snl}} = \sum_{n=M+1}^{\infty} \frac{\mu(n)}{n} \log(\zeta(sn))$$

Taking in both sides the absolute value we get from Proposition 1

$$\left| \sum_{k=1}^{\infty} \frac{1}{p_k^s} - \sum_{n=1}^{M} \frac{\mu(n)}{n} \sum_{k=1}^{\infty} \sum_{l=1}^{\infty} \frac{1}{l} \frac{1}{p_k^{snl}} \right| =$$

$$= \left| \sum_{n=M+1}^{\infty} \frac{\mu(n)}{n} \log(\zeta(sn)) \right|$$

$$\leq C(f) \frac{2^s (2^s + 1)(s + 1)}{(2^s - 1)^3} \frac{M^2}{(sM + s - 1)2^{sM}}$$

$$\leq C(f) \frac{M}{2^{s(M+1)}}$$

Let now $\dfrac{g^{(a)}(0)}{a!}$, $a = 1, 2, 3,\ldots$ are bounded above. We have

$$\dfrac{g^{(a)}(0)}{a!}\left(\sum_{k=1}^{\infty}\dfrac{1}{p_k^{sa}} - \sum_{n=1}^{M}\dfrac{\mu(n)}{n}\sum_{k=1}^{\infty}\sum_{l=1}^{\infty}\dfrac{1}{l}\dfrac{1}{p_k^{snla}}\right) = O\left(\dfrac{M}{2^{sa(M+1)}}\dfrac{g^{(a)}(0)}{a!}\right)$$

We sum the above estimation with respect to $a$, since $\dfrac{g^{(a)}(0)}{a!}$ are bounded above we get $\sum_{k=1}^{\infty}g\left(\dfrac{1}{p_k^s}\right) - \sum_{n=1}^{M}\dfrac{\mu(n)}{n}\sum_{k=1}^{\infty}\sum_{l=1}^{\infty}\dfrac{1}{l}g\left(\dfrac{1}{p_k^{snl}}\right) = O\left(M \cdot g\left(\dfrac{1}{2^{s(M+1)}}\right)\right)$ and the result follows replacing $M$ with $x$. □